\documentclass[12pt]{amsart}
\usepackage{color}

\newcommand{\reals}{{\mathbb R}}
\newcommand{\nats}{{\mathbb N}}
\newcommand{\prob}{{\mathbb P}}
\newcommand{\expe}{{\mathbb E}}
\DeclareMathOperator{\inte}{int}
\newtheorem{lemma}{Lemma}
\newtheorem{theorem}{Theorem}
\newtheorem{proposition}{Proposition}
\newtheorem{corollary}{Corollary}
\theoremstyle{definition}
\newtheorem{example}{Example}
\newtheorem{rem}{Remark}

\begin{document}
\begin{center}
\Large Bayesian Posteriors For Arbitrarily Rare Events
\end{center}

\vspace*{10mm}

\begin{center}
Drew Fudenberg$^a$, Kevin He$^b$, and Lorens A.\ Imhof$^c$
\end{center}

\vspace*{10mm}

\noindent $^a$Department of Economics,
Massachusetts Institute of Technology, Cambridge, MA 02139, USA,
E-mail: drew.fudenberg@gmail.com 

\vspace*{5mm}

\noindent $^b$Department of Economics, Harvard University, 
Cambridge, MA 02138, USA, 
\\E-mail: hesichao@gmail.com

\vspace*{5mm}

\noindent $^c$Department of Statistics and Hausdorff Center for Mathematics,
Bonn University, 53113 Bonn, Germany, E-mail: limhof@uni-bonn.de

\vspace*{10mm}


\vspace*{10mm}

Abstract.
We study how much data a Bayesian observer needs to correctly infer the relative likelihoods of two events when both events are arbitrarily rare. Each period, either a blue die or a red die is tossed. The two dice
land on side $1$ with unknown probabilities $p_1$ and $q_1$, which can 
be arbitrarily low. Given a data-generating process where $p_1\ge c q_1$, 
we are interested in how much data is required to guarantee that with 
high probability the observer's Bayesian posterior mean for $p_1$ exceeds
$(1-\delta)c$ times that for $q_1$. If the prior densities for the two 
dice are positive on the interior of the parameter space and behave like 
power functions at the boundary, then for every $\epsilon>0,$ there exists 
a finite $N$ so that the observer obtains such an inference after $n$ 
periods with probability at least $1-\epsilon$ whenever $np_1\ge N$. The 
condition on $n$ and $p_1$ is the best possible. The result can fail 
if one of the prior densities converges to zero exponentially fast 
at the boundary.

\vspace*{10mm}

Significance Statement.
Many decision problems in contexts ranging from drug safety tests to game-theoretic learning models require Bayesian comparisons between the likelihoods of two events. When both events are arbitrarily rare, a large data set is needed to reach the correct decision with high probability. The best result in previous work requires the data size to grow so quickly with rarity that the expectation of the number of observations of the rare event explodes. We show for a large class of priors that it is enough that this expectation exceeds a prior-dependent constant. However, without some restrictions on the prior the result fails, and our 
condition on the data size is the weakest possible.

\section{Introduction}
Suppose a physician is deciding between 
a routine surgery versus a newly approved drug for her patient.
Either treatment can, in rare cases, lead to a life-threatening 
complication. She adopts a Bayesian approach to estimate the respective 
probability of complication, as is common among practitioners in medicine 
when dealing with rare events, 
see, for example, 
U.S. Food and Drug Administration (2000) and Thompson (2014) on 
the ``zero-numerator problem.'' She reads the medical literature to learn 
about $n$ patient outcomes associated with the two treatments and chooses the 
new drug if and only if her posterior mean regarding the probability of 
complication due to the drug is lower than $(1-\delta)$ times that of the 
surgery. As the true probability of complication becomes small 
for both treatments, how quickly does $n$ need to increase to ensure that
the physician will correctly choose surgery with probability at least
$1-\epsilon$ when surgery is in fact the safer option?

Phrased more generally, we study how much data is required for the Bayesian 
posterior means on two probabilities to respect an inequality between them 
in the data-generating process, where these true probabilities may be 
arbitrarily small. Each period, one of two dice, blue or red, is chosen to be 
tossed. The choices can be deterministic or random, but have to be independent 
of past outcomes. The blue and red dice land on side $k$ with unknown 
probabilities $p_{k}$ and $q_{k}$, and the outcomes of the tosses 
are independent of past outcomes. Say that the posterior beliefs 
of a Bayesian observer satisfy \emph{$(c,\delta)$-monotonicity}
for side $\bar{k}$ if his posterior mean for $p_{\bar{k}}$ exceeds $(1-\delta) c $
times that for $q_{\bar{k}}$ whenever the true probabilities are such that 
$p_{\bar{k}}\ge c q_{\bar{k}}$. 
We assume the prior densities are continuous and positive on the interior of
the probability simplex, and behave like power functions at the boundary.
Then we show that, under a mild condition on the frequencies of the chosen colors, 
for every $\epsilon>0$, there exists a finite $N$
so that the observer holds a $(c,\delta)$-monotonic belief after $n$
periods with probability at least $1-\epsilon$ whenever $np_{\bar{k}}\geq N$.
This condition means that the expected number of times the blue
die lands on side $\bar k$ must exceed a constant that is independent
of the true parameter.
Examples show that the sample size condition is the best possible, and that the 
result can fail if one of the prior densities converges to zero exponentially 
fast at the boundary. 
A crucial aspect of our problem is the behavior of estimates
when the true parameter value approaches the boundary of the parameter
space, a situation that is rarely studied in a Bayesian context. 

Suppose that in every period, the blue die is chosen with the same probability 
and that outcome $\bar k$ is more likely under the blue die than under the red one.  
Then, under our conditions, an observer who sees outcome $\bar k$ but not the die color is 
very likely to assign posterior odds ratio to blue versus red that is not much below the 
prior odds ratio. That is, the observer is unlikely to update her beliefs in the wrong 
direction. Fudenberg and He (2017) use this corollary to provide a learning-based foundation 
for equilibrium refinements in signalling games.

The best related result known so far is a consequence of the uniform consistency 
result of Diaconis and Freedman (1990). Their result leads to the desired conclusion 
only under the stronger condition that the sample size is so large that the expected 
number of times the blue die lands on side $\bar k$ exceeds a threshold proportional 
to $1/p_{\bar{k}}$. That is, the threshold obtained from their result explodes as 
$p_{\bar{k}}$ approaches zero.

Our improvement of the sample size condition is made possible by
a new pair of inequalities that relate the 
Bayes estimates to observed frequencies. Like the bounds of Diaconis 
and Freedman (1990), the inequalities apply to all sample sequences 
without exceptional null sets and they do not involve true parameter values.
Our result is related to a recent result of Bochkina and Green (2014) which shows
that, under some conditions, the posterior distribution converges 
faster when the true parameter is on the boundary.
Our result is also related to Dudley and Haughton (2002) who
consider a half-space not containing the maximum likelihood estimate
of the true parameter, and study how quickly the posterior probability 
assigned to the half-space converges to zero.

\section{Bayes estimates for multinomial probabilities}
\label{sectiononemultinomial}

We first consider the simpler problem of estimating for a single 
$K$-sided die the probabilities of landing on the various sides.
Suppose the die is tossed independently $n$ times. Let $X_k^n$ denote
the number of times the die lands on side $k$. Then
$X^n=(X_1^n,\dots,X_K^n)$ has
a multinomial distribution with parameter $n\in\nats$ and unknown parameter 
$p=(p_1,\dots,p_K)\in\Delta$, where 
$\nats$ is the set of positive integers and
$\Delta=  \{p\in[0,1]^K : p_1+\dots +p_K=1\}$.
Let ${\nats}_0 = \nats \cup \{0\}$.
Let $\pi$ be a prior density on $\Delta$ with respect to 
the Lebesgue measure $\lambda $ on $\Delta$,
normalized by $\lambda(\Delta)=1/(K-1)!$.
Let $\pi(\cdot | X^n)$ be the posterior density after observing $X^n$.

Motivated by applications where some of the $p_k$ can be arbitrarily small, 
we are interested in whether the {\em relative} error of the Bayes estimator
$\hat p_k(X^n)=\int p_k \pi(p|X^n)\, d\lambda(p)$ is small with probability 
close to $1$, uniformly on large subsets of $\Delta$. Specifically, given 
$k\in\{1,\dots, K\}$ and $\epsilon>0$, we seek conditions on $n$ and $p$  
and the prior, so that 
\begin{equation}\label{diceineq}
\prob_p(|\hat p_k(X^n) - p_k| < p_k\epsilon ) \geq 1-\epsilon.
\end{equation}
A subscript on $\prob$ or $\expe$ indicates the parameter value
under which the probability or expectation is to be taken.


For a wide class of priors, we show in Theorem~\ref{thm1} that there is a constant 
$N$ that is independent of the unknown parameter so that (\ref{diceineq}) holds 
whenever $\expe_p(X_k^n)\geq N$. Denote the interior of $\Delta$ by $\inte\Delta$.

\medskip

{\bf Condition $\boldsymbol {\mathcal P}$.} We say that a density $\pi$ on $\Delta$ satisfies Condition
${\mathcal P}(\alpha)$, where $\alpha=(\alpha_1,\dots,\alpha_K)\in(0,\infty)^K$, if 
\[
\frac{\pi(p)}{\prod_{k=1}^K p_k^{\alpha_k-1}}
\]
is uniformly continuous and bounded away from zero on $\inte\Delta$.
We say that $\pi$ satisfies Condition $\mathcal P$ if there exists $\alpha\in(0,\infty)^K$ so that $\pi$
satisfies Condition ${\mathcal P}(\alpha)$.

\medskip

For example, if $K=2$, then $\pi$ satisfies Condition ${\mathcal P}(\alpha)$ if and only if $\pi$ is
positive and continuous on $\inte\Delta$ and the limit $\lim_{p_k\to 0}\pi(p)/p_k^{\alpha_k-1}$
exists and is positive for $k=1,2$. For every $K\geq 2$, every Dirichlet distribution
has a density that satisfies Condition $\mathcal P$.
Note that Condition ${\mathcal P}$ does not require that the density
is bounded away from zero and infinity at the boundary.
The present assumption on the behavior at the boundary is similar
to Assumption P of Bochkina and Green (2014).

\begin{theorem}\label{thm1}
Suppose $\pi$ satisfies Condition ${\mathcal P}$.
Then for every $\epsilon>0$, there exists $N\in\nats$ so that
\begin{equation}\label{thm1ineq}
\prob_p(\left|\hat p_k(X^n) - p_k \right| \geq p_k \epsilon) \leq \epsilon
\end{equation}
if $np_k \geq N$.
\end{theorem}

The proofs of the results in this section are given in the Supporting Information.

The proof of Theorem~\ref{thm1} uses bounds on the posterior means given in 
Proposition~\ref{prop1} below. These bounds imply that there is an $N\in\nats$ so that if 
$np_k \geq N$ and the maximum likelihood estimator $\frac 1n X^n_k$ is close to $p_k$, 
then $|\hat p_k(X^n) -p_k |< p_k \epsilon$. It follows from Chernoff's inequality that 
the probability that $\frac 1n X^n_k$ is not close to $p_k$ is at most $\epsilon$.

Inequality (\ref{thm1ineq}) shows a higher accuracy of the Bayes estimator $\hat p_k(X^n)$ 
when the true parameter $p_k$ approaches $0$. To explain this fact in a special case suppose 
that $K=2$ and the prior is the uniform distribution. Then $\hat p_k(X^n)=(X_k^n+1)/(n+2)$ 
and the mean squared error of $\hat p_k(X^n)$ is $[np_k(1-p_k)+(1-2p_k)^2]/(n+2)^2$, which 
converges to $0$ like $\frac 1n$ when $p_k\in(0,1)$ is fixed, and like $\frac 1{n^2}$ when 
$p_k = \frac{1}{n}$. Moreover, by Markov's inequality, the probability in (\ref{thm1ineq})
is less than $(np_k+1)/(n^2p_k^2\epsilon^2)$, so that in this case we can choose 
$N= 2/\epsilon^3$. In general, we do not have an explicit expression for the threshold $N$,
but in Remark~\ref{rem2} we discuss the properties of the prior that have an impact
on the $N$ we construct in the proof.

Condition $\mathcal P$ allows the prior density to converge to zero at the boundary of 
$\Delta$ like a power function with an arbitrarily large exponent. The following example 
shows that the conclusion of Theorem~\ref{thm1} fails to hold for a prior density
that converges to $0$ exponentially fast. 

\begin{example}\label{ex1}
Let $K=2$,  $\pi(p)\propto e^{-1/p_1}$, and $\delta>0$. Then for every $N\in\nats$, there
exist $p\in\Delta$ and $n\in\nats$ with $n^{\frac 12+\delta} p_1 \geq N$  so that
\[
\prob_p( |\hat p_1(X^n) - p_1| > p_1) =1.
\]
\end{example}

The idea behind this example is that the prior assigns very little
mass near the boundary point where $p_1=0$, so if the true parameter $p_1$ is small, 
the observer needs a tremendous amount of data to 
be convinced that $p_1$ is in fact small. The prior density in our example converges to 
$0$ at an exponential rate as $p_1\to 0$, and it turns out that the amount of data 
needed in order that $\hat p_1(X^n)/p_1$ is close to $1$ grows quadratically in $1/p_1$.  
For every fixed $N\in\nats$ and $\delta>0$, the pairs $(n,p_1)$ satisfying the relation 
$n^{\frac 12+\delta} p_1 =N$ involve a sub-quadratic growth rate of $n$ with respect to 
$1/p_1$. So we can always pick a small enough $p_1$ such that the corresponding data size 
$n$ is insufficient.

The next example shows that the sample size condition of
Theorem~\ref{thm1}, $np_k\geq N$, cannot be replaced by a weaker condition
of the form $\zeta(n) p_k \geq N$ for some function $\zeta$ with
$\limsup_{n\to\infty}\zeta(n)/n=\infty$. Put differently, the set of
$p$ for which (\ref{thm1ineq}) can be proved cannot be enlarged to a set 
of the form $\{p : p_k \geq \phi_\epsilon(n)\}$ with $\phi_\epsilon(n)=o(1/n)$.

\begin{example}\label{ex2}
Suppose $\pi$ satisfies Condition $\mathcal P$. 
Let $\zeta: \nats\to (0,\infty)$ be so that $\limsup_{n\to\infty}\zeta(n)/n=\infty$.
Then for every $N\in\nats$, there exist $p\in\Delta$ and $n\in\nats$ with 
$\zeta(n)p_1\geq N$ so that
\[
\prob_p(|\hat p_1(X^n) - p_1 | > p_1) =1.
\]
\end{example}

The following proposition gives fairly sharp bounds on the posterior means
under the assumption that the prior density satisfies Condition $\mathcal P$.
The result is purely deterministic and applies to all possible sample sequences.
The bounds are of interest in their own right, and also play
a crucial role in the proofs of Theorem~\ref{thm1} and \ref{thm2}.

\begin{proposition}\label{prop1} 
Suppose $\pi$ satisfies Condition ${\mathcal P}(\alpha)$.
Then for every $\epsilon>0$, there exists a constant $\gamma>0$ such that
\begin{equation}\label{prop1ineq}
(1-\epsilon)\frac {n_k+\alpha_k} {n+ \gamma}
\leq
\frac{\int p_k \left(\prod_{i=1}^K p_i^{n_i}\right) \pi(p) \, d\lambda(p)}
{\int\left(\prod_{i=1}^K p_i^{n_i}\right)\pi(p)\, d\lambda(p)}
\leq 
(1+\epsilon) \frac{n_k+\gamma} {n+\gamma} 
\end{equation}
for $k=1,\dots,K$ and all $n, n_1,\dots,n_K\in\nats_0$ with $\sum_{i=1}^Kn_i=n$. 
\end{proposition}

\begin{rem}\label{rem1}
If $\pi$ is the density of a Dirichlet distribution with parameter $\alpha\in(0,\infty)^K$, 
then the inequalities in (\ref{prop1ineq}) hold with $\epsilon=0$ and 
$\gamma=\sum_{k=1}^K\alpha_k$, and the inequality on the left-hand side is an equality.
If $\pi$ is the density of a mixture of Dirichlet distributions and the support of the 
mixing distribution is included in the interval $[a,A]^K$, $0\leq a \leq A <\infty$, then 
for all $k$ and $n_1,\dots,n_K$ with $\sum_{i=1}^Kn_i=n$,
\begin{equation}\label{mixineq}
\frac{n_k +a}{n+ K A}
\leq\frac{\int p_k \left(\prod_{i=1}^K p_i^{n_i}\right) \pi(p)\, d\lambda(p)}
{\int \left(\prod_{i=1}^K p_i^{n_i} \right) \pi(p)\, d\lambda(p)}
\leq\frac {n_k + A}{n+ Ka}.
\end{equation}
 
The proofs of our main results, Theorems~\ref{thm1} and \ref{thm2}, apply to all 
priors whose densities satisfy inequalities (\ref{prop1ineq}) or (\ref{mixineq}).
In particular, the conclusions of these theorems and of their corollaries hold if the prior 
distribution is a mixture of Dirichlet distributions and the support of the mixing 
distribution is bounded.
\end{rem}

\begin{rem}\label{rem2}
Condition $\mathcal P(\alpha)$ implies that the function $\pi(p)/\prod_{k=1}^Kp_k^{\alpha_k-1}$,
$p\in\inte\Delta$, can be extended to a continuous function $\tilde\pi(p)$ on $\Delta$.
The proof of Proposition~\ref{prop1} relies on the fact that $\tilde\pi$ can be uniformly 
approximated by Bernstein polynomials. 
An inspection of the proof shows that the constant $\gamma$ in (\ref{prop1ineq})
can be taken to be $m +\sum_{k=1}^K\alpha_k$, where $m$ is so large that 
$h_m$, the $m$-th degree Bernstein polynomial of $\tilde\pi$, satisfies
\[
\max \{|h_m(p) -\tilde\pi(p)| :p\in\Delta\} \leq 
\frac{\min\{\tilde\pi(p) : p\in\Delta\} }{1+2\epsilon^{-1}}.
\]
Hence, in addition to a small value of $\epsilon$, the following
properties of the density $\pi$ result in a large value of $\gamma$:
(i) if $\sum_{k=1}^K\alpha_k$ is large, (ii) if $\pi$ is a ``rough'' function 
so that $\tilde\pi$ is hard to approximate and $m$ needs to be large, 
and (iii) if $\tilde\pi$ is close to $0$ somewhere. The threshold $N$
in Theorem~\ref{thm1} depends on the prior through the constant $\gamma$
from Proposition~\ref{prop1} and the properties of $\pi$ just described will
also lead to a large value of $N$. 

In particular, $N\to\infty$ if $\sum_{k=1}^K\alpha_k\to\infty$. 
For example,
consider a sequence of priors $\pi^{(j)}$ for $K=2$ where $\pi^{(j)}$
is the density of the Dirichlet distribution with parameter $(j,1)$, 
so that $\pi^{(j)}$ satisfies Condition ${\mathcal P}(\alpha)$ with $\alpha_1 = j$. 
As $j\to \infty$, $\pi^{(j)}$ converges faster and faster to $0$ as $p_1 \to 0$, 
though never as fast as in Example~\ref{ex1}, where no finite $N$ can satisfy 
the conclusion of Theorem~\ref{thm1}. If $n=4j$ and $p_1=\frac 1{12}$, then under $\pi^{(j)}$,
$\hat p_1(X^n) = (X_1^n+j)/(n+j+1) \geq 2p_1$, so for every $\epsilon\in(0,1)$, the
probability in Theorem~\ref{thm1} is $1$. Thus, the smallest $N$ for which the conclusion 
holds must exceed $4j\times \frac 1{12} = \frac j3$.
\end{rem}

\begin{rem}
Using results on the degree of approximation by Bernstein polynomials,
one may compute explicit values for the constants $\gamma$ in Proposition~\ref{prop1}
and $N$ in Theorem~\ref{thm1}. Details are given in Remarks $3'$
and $3''$ in Supporting Information.
\end{rem}

\begin{rem}\label{rem4}
Suppose $K>2$ and the statistician is interested in only one of the probabilities $p_k$, 
say $p_{\bar k}$. Then, instead of using $\hat p_{\bar k}(X^n)$, he may first reduce the 
original $(K-1)$-dimensional estimation problem to the problem of estimating the 
one-dimensional parameter 
$(p_{\bar k} ,\sum_{k\neq {\bar k}} p_k)$ of the Dirichlet 
distribution of $(X_{\bar k}^n ,\sum_{k\neq {\bar k}} X_k^n)$. He will then 
distinguish only whether or not the die lands on side $\bar k$ and will use the induced 
one-dimensional prior distribution for the parameter of interest.
If the original prior is a Dirichlet distribution on $\Delta$, both approaches lead 
to the same Bayes estimators for $p_{\bar k}$, but in general, they do not.
Proposition~\ref{prop2} in the Supporting Information shows that whenever 
the original density $\pi$ satisfies condition $\mathcal P$, then the induced 
density satisfies condition $\mathcal P$ as well.
However, it may happen that the induced density satisfies Condition 
$\mathcal P$ even though the original density does not.
For example, if $K=3$ and $\pi(p)\propto e^{-1/p_1} + p_2$, then $\pi$ does not satisfy 
Condition $\mathcal P$, but for each $\bar k=1,2,3$, the induced density does. 
\end{rem}

\section{Comparison of two multinomial distributions}
\label{sectiontwomultinomials}

Here we consider two dice, blue and red, each with $K\geq 2$ sides. In every period,
a die is chosen. We first consider the case where the choice is deterministic and 
fixed in advance. We will later allow the choice to be random.
The chosen die is tossed and lands on the $k$-th side according to the unknown
probability distributions $p=(p_1,\dots,p_K)$ and $q=(q_1,\dots,q_K)$
for the blue and the red die, respectively. 
The outcome of the toss is independent of past outcomes.
The parameter space of the
problem is $\Delta^2$. 
The observer's prior is represented by a product density 
$\pi(p)\varrho(q)$ over $\Delta^2$, that is, he regards the
parameters $p$ and $q$ as realizations of independent random vectors. 

Let $X^n$ be a random vector that describes the outcomes, i.e., colors and sides, of the 
first $n$ tosses. Let $b_n$ denote the number of times the blue die is tossed in the first 
$n$ periods. Let $\pi(\cdot | X^n)$ and $\varrho(\cdot| X^n)$ be the posterior densities
for the blue and the red die after observing $X^n$. Let 
$\hat p_k(X^n)=\int p_k \pi(p|X^n)\, d\lambda(p)$ and 
$\hat q_k(X^n)=\int q_k \varrho(q|X^n)\, d\lambda(q)$.
The product form of the prior density ensures that the marginal posterior distribution for
either die is completely determined by the observations on that die and the marginal prior
for that die.

We study the following problem. Fix a side $\bar k\in\{1,\dots,K\}$ and a constant 
$c\in(0,\infty)$. Consider a family of environments, each characterized by a data-generating 
parameter vector $\vartheta=(p,q)\in\Delta^2$ and an observation length $n$. In each 
environment, we have $p_{\bar k}\geq c q_{\bar k}$, and we are interested in whether the 
Bayes estimators reflect this inequality. In general, one cannot expect that the probability 
that $\hat p_{\bar k}(X^n) \geq c \hat q_{\bar k}(X^n)$ is much higher than $\frac 12$ when 
$p_{\bar k}=c q_{\bar k}$. We therefore ask whether in all of the environments, the observer 
has a high probability that $\hat p_{\bar k}(X^n) \geq c(1-\delta) \hat q_{\bar k}(X^n)$ for 
a given constant $\delta\in(0,1)$.

Clearly, as $p_{\bar k}$ approaches $0$, we will need a larger observation length $n$ for
the data to overwhelm the prior. But how fast must $n$ grow relative to $p_{\bar k}$? 
Applying the uniform consistency result of Diaconis and Freedman to each Bayes estimator 
separately leads to the condition that $n$ must be so large that the expected number of 
times the blue die lands on side $\bar k$, that is, $b_n p_{\bar k}$, exceeds a threshold 
that explodes when $p_{\bar k}$ approaches zero. The following theorem shows that there is 
a threshold that is independent of $p$, provided the prior densities satisfy Condition 
$\mathcal P$.

\begin{theorem}\label{thm2}
Suppose that $\pi$ and $\varrho$ satisfy Condition $\mathcal P$.
Let $\bar k\in\{1,\dots,K\}$, $c\in(0,\infty)$, and $\epsilon,\delta,\eta\in (0,1)$. Then 
there exists $N\in\nats$ so that for every deterministic sequence of 
choices of the dice to be tossed,
\begin{equation}\label{thm2ineq}
\prob_\vartheta (\hat p_{\bar k} (X^n) \geq c (1-\delta) \hat q_{\bar k}(X^n)) \geq 1-\epsilon
\end{equation}
for all $\vartheta=(p,q)\in\Delta^2$ with $p_{\bar k} \geq c q_{\bar k} $ and 
all $n\in\nats$ with $b_n p_{\bar k} \geq N$ and $b_n/n \leq  1-\eta $.
\end{theorem}

We prove Theorem~\ref{thm2} in the next section.

Note that the only constraints on the sample size here are that the product of $b_n$ with 
$p_{\bar{k}}$ be sufficiently large and the proportion of periods in which the red die is 
chosen be not too small. However, $p_{\bar{k}}$ and $q_{\bar{k}}$ can be arbitrarily small. 
This is useful in analyzing situations where the data-generating process contains rare events.

In the language of hypothesis testing, Theorem~\ref{thm2} says that under the stated condition 
on the prior, the test that rejects the null hypothesis $p_{\bar k}\geq c q_{\bar k}$ if and 
only if $\hat p_{\bar k} (X^n) < c (1-\delta) \hat q_{\bar k}(X^n)$ has a type I error 
probability of at most $\epsilon$ provided $p_{\bar k}\geq N/ b_n$ (and $b_n/n < 1-\eta $).
For every $n$, the bound on the error probability holds uniformly on the specified parameter 
set. Note that such a bound cannot be obtained for a test that rejects the hypothesis whenever 
$\hat p_{\bar k} (X^n) < c \hat q_{\bar k}(X^n)$.

We now turn to the case where the dice are randomly chosen. 
The probability of choosing the blue die need not be constant over time but
must not depend on the unknown parameter $\vartheta$.
Let the random variable 
$B_n$ denote the number of times the blue die is tossed in the first $n$ periods.

\begin{corollary}\label{cor1}
Suppose that $\pi$ and $\varrho$ satisfy Condition $\mathcal P$.
Let $\bar k\in\{1,\dots,K\}$, $c\in(0,\infty)$, and $\epsilon,\delta\in (0,1)$. 
Suppose that in every period, the die to be tossed is chosen at random, independent
of the past, 
and that 
\begin{equation}\label{EBncondition}
\liminf_{n\to\infty}\frac {\expe(B_n)}n >0,\qquad \limsup_{n\to\infty}\frac {\expe(B_n)}n <1.
\end{equation}
Then there exists $N\in\nats$ so that
\begin{equation}\label{cor1ineq}
\prob_\vartheta (\hat p_{\bar k} (X^n) \geq c (1-\delta) \hat q_{\bar k}(X^n)) \geq 1-\epsilon
\end{equation}
for all $\vartheta=(p,q)\in\Delta^2$ with $p_{\bar k} \geq c q_{\bar k} $ and 
all $n\in\nats$ with $n p_{\bar k} \geq N$.
\end{corollary}

The proof of Corollary~\ref{cor1} is given at the end of the next section.

In the decision problem described in the first paragraph of the introduction, Theorem~\ref{thm2} 
and Corollary~\ref{cor1} ensure that whenever surgery is the safer option, the probability that 
the physician actually chooses surgery is at least $1-\epsilon$ unless the probability of 
complication due to the drug is smaller than $N/n$. Except for this last condition, the bound 
$1-\epsilon$ holds uniformly over all possible parameters.

In the rest of this section we assume that in every period the blue die is chosen at random with 
the same probability $\mu_B$. The value of $\mu_B$ need not be known, we only assume that 
$0<\mu_B<1$, so that condition (\ref{EBncondition}) is met.

The following example shows that the conditions on the prior densities cannot be omitted from 
Corollary~\ref{cor1}. 

\begin{example}\label{ex3} Suppose $K=2$ and $0<\mu_B<1$. Suppose $\pi$ 
satisfies Condition $\mathcal P$ and $\varrho(q)\propto e^{-1/q_1}$. Let 
$c>0$. Then for every $N\in\nats$, there exist $\vartheta=(p,q)\in\Delta^2$
with $p_1\geq c q_1$ and $n\in\nats$ with $np_1\geq N$ so that 
\[
\prob_{\vartheta}\left(\hat p_1(X^n) < \frac c2 \hat q_1(X^n)\right) > \frac{1}{2}.
\]
\end{example}

The next example shows that the sample size condition of Corollary~\ref{cor1},
$np_{\bar{k}}\ge N$, is the best possible for small $p_{\bar{k}}$. It cannot be 
replaced by a weaker condition of the form $n\zeta(p_{\bar{k}})\ge N$ for some
function $\zeta$ with $\lim_{t\to0+}\zeta(t)/t=\infty$. In particular, taking $\zeta$ 
to be a constant function shows that there does not exist $N\in\nats$ so that $n\geq N$ 
implies that (\ref{cor1ineq}) holds uniformly for all $\vartheta$ with $p_k \geq cq_k$.

\begin{example}\label{ex4}
Suppose that $0<\mu_B<1$ and that $\pi$  and $\varrho$ satisfy Condition 
$\mathcal P$. Let $c>0$. Let $\zeta$ be a nonnegative function on $[0,1]$ 
with $\lim_{t\to 0+} \zeta(t)/t=\infty$.
Then there exists $\epsilon_0>0$ so that for every $N\in\nats$, there exist 
$\vartheta=(p,q)\in\Delta^2$ with $p_1\geq c q_1$ and $n\in\nats$ with $n\zeta(p_1) \geq N$
so that
\[
\prob_{\vartheta}\left( \hat p_1(X^n) < \frac c2 \hat q_1(X^n) \right) > \epsilon_0.
\]
\end{example}

Examples~\ref{ex3} and \ref{ex4} are proved in the Supporting Information.

Suppose that after data $X^{n}$ the observer were told that the next
outcome was $\bar{k}$ but not which die was used. Then Bayes' rule
implies the posterior odds ratio for ``blue'' relative to ``red'' is 
\[
\frac{\mu_{B}\int p_{\bar{k}}\pi(p|X^{n})\, d\lambda(p)}
{\mu_{R}\int q_{\bar{k}}\varrho(q|X^{n})\, d\lambda(q)},
\]
where $\mu_R=1-\mu_B$.

\begin{corollary}\label{cor2}
Suppose that $\pi$ and $\varrho$ satisfy Condition $\mathcal P$. Then there 
exists $N\in\mathbb{N}$ such that whenever $p_{\bar{k}}\geq q_{\bar{k}}$ and 
$np_{\bar{k}}\geq N$, there is probability at least $1-\epsilon$ that the 
posterior odds ratio of ``blue'' relative to ``red'' exceeds 
$(1-\epsilon)\cdot\frac{\mu_{B}}{\mu_{R}}$ when the $(n+1)$-th die lands on 
side $\bar{k}$. 
\end{corollary}

This corollary is used in Fudenberg and He (2017), who provide a learning-based foundation for equilibrium refinements in signalling games. They consider a sequence of learning environments, each containing populations of ``blue'' senders, ``red'' senders, and receivers. Senders are randomly matched with receivers each period and communicate using one of $K$ messages. There is some special message $\bar k$, whose probability of being sent by blue senders always exceeds the probability of being sent by red senders in each environment. Suppose the common prior of the receivers satisfies Condition $\mathcal P$ and, in every environment, there are enough periods that the expected total observations of blue sender playing $\bar k$ exceeds a constant. Then at the end of every environment, by this corollary all but $\epsilon$ fraction of the receivers will assign a posterior odds ratio for the color of the sender not much less than the prior odds ratio of red versus blue, if they were to observe another instance of $\bar k$ sent by an unknown sender, regardless of how rarely the message $\bar k$ is observed.  A leading case of receiver prior satisfying Condition 
$\mathcal P$ is fictitious play, the most commonly used model of learning in games, which corresponds to Bayesian updating from a Dirichlet prior, but our corollary shows that the Dirichlet restriction can be substantially relaxed.

\section{Proofs of Theorem~\ref{thm2} and Corollary~\ref{cor1}}\label{proofofthm2}
We begin with two auxiliary results needed in the proof of Theorem~\ref{thm2}.
Lemma~\ref{la3} is a large deviation estimate which gives a bound on the probability 
that the frequency of side $\bar k$ in the tosses of the red die exceeds an affine 
function of the frequency of side $\bar k$ in the tosses of the blue die. 
Lemma~\ref{la4} implies that, with probability close to $1$, the number of times the 
blue die lands on side $\bar k$ exceeds a given number when $b_np_{\bar k}$ is 
sufficiently large. The proofs of Lemmas~\ref{la3} and ~\ref{la4} are in the
Supporting Information.

\begin{lemma}\label{la3}
Let $S_n$ be a binomial random variable with parameters $n$ and $p$, and let $T_m$ be a binomial 
random variable with parameters $m$ and $q$. Let $0<c'<c$ and $d>0$.
Suppose $S_n$ and $T_m$ are independent, and $p\geq c q$.
Then
\[
\prob\left(\frac {T_m}m  \geq  \frac 1{c'} \frac{S_n}n + \frac d{n\wedge m} \right) 
\leq \left( \frac {c'}c \right)^{c'd /(c' +1)}.
\]
\end{lemma}

\begin{lemma}\label{la4}
Let $M<\infty$ and $\epsilon>0$. 
Then there exists $N\in\nats$ so that if $S_n$ is a binomial random 
variable with parameters $n$ and $p$ and $n p \geq N$, then
\[
\prob_p(S_n \leq  M) \leq \epsilon.
\]
\end{lemma}

{\em Proof of Theorem~\ref{thm2}.}
Let $r_n=n-b_n$ be the number of times the red die is tossed in the first $n$ periods. 
Let $Y_n$ and $Z_n$ be the respective number of times the blue 
and the red die land on side $\bar k$. Choose $\beta>0$ and $c'\in(0,c)$ so that
\begin{equation}\label{beta-c'-choice}
\frac{1-\beta}{(1+\beta)(1-\delta)} > \frac c{c'}+\delta .
\end{equation}
By Proposition~\ref{prop1}, there exists $\gamma>0$ so that for every $n\in\nats$,
\begin{equation}\label{bounds}
\hat p_{\bar k}(X^n) \geq \phi(b_n ,Y_n),\qquad \hat q_{\bar k}(X^n) \leq \psi(r_n, Z_n),
\end{equation}
where
\[ 
\phi(b,y)=(1-\beta)\frac y{b+\gamma}, 
\qquad
\psi(r,z)=(1+\beta)\frac {z+\gamma}{r}.
\]
Let $d>0$ be so that the bound in Lemma~\ref{la3} satisfies 
$\left(c'/c\right)^{c'd/(c' + 1)} \leq \frac \epsilon 2$.

We now show that for all $b,r\in\nats$, $y=0,\dots,b$, and $z=0,\dots, r$, the inequalities
\begin{equation}\label{ineq1}
\frac zr < \frac{1}{c'}\frac yb + \frac{d}{b\wedge r},
\quad
\frac {2c\gamma}{c' \delta} < b< \frac{r}{\eta},
\quad 
y>M:=\frac{3c(d+\gamma)}{\delta\eta}
\end{equation}
imply that 
\begin{equation}\label{ineq2}
\phi(b,y) > c (1-\delta) \psi(r,z).
\end{equation}
It follows from the first and the third inequality in (\ref{ineq1}) that
\begin{align*}
\psi(r,z) & < \psi\left(r, \frac {ry}{c'b} + \frac{rd}{b\wedge r}\right)\\
&=(1+\beta)\left(\frac y{c'b}+ \frac{d}{b\wedge r} + \frac{\gamma}{r} \right)\\
& \leq (1+\beta) \left( \frac y{c'b} + \frac{\delta M}{3bc}\right).
\end{align*}
Applying this result, inequality (\ref{beta-c'-choice}), twice the second and finally the 
fourth inequality in (\ref{ineq1}) we get
\begin{align*}
& \frac{\phi(b,y) -c(1-\delta)\psi(r,z)}{(1-\delta)(1+\beta)}\\
& \qquad\quad >
\frac{y}{b+\gamma} \left( \frac{1-\beta}{(1-\delta)(1+\beta)} - \frac c{c'} - \frac{c \gamma}{c'b}  \right)  - \frac{\delta M}{3b}\\
&\qquad \quad\geq \frac{y}{b+\gamma} \left( \delta - \frac\delta 2 \right)  - \frac{\delta M}{3b}\\
&\qquad \quad\geq \frac{2}{3b} \frac \delta 2 M - \frac{\delta M}{3b}=0,
\end{align*}
proving (\ref{ineq2}).

Let ${\mathcal N} = \{n\in\nats : b_n/n \leq 1-\eta\}$, $N_1= \lceil 2c\gamma / (c'\delta)\rceil$, 
and for every $n\in{\mathcal N}$ with $b_n\geq N_1$ define the events
\[
F_n = \left\{\frac {Z_n}{r_n} < \frac 1 {c'} \frac {Y_n}{b_n} + \frac d{b_n \wedge r_n} \right\},
\qquad
G_n = \{Y_n > M\}.
\]
For all $n\in{\mathcal N}$, $b_n < r_n/\eta$. Thus, if $n\in{\mathcal N}$ and $b_n\geq N_1$, 
the implication  $(\ref{ineq1}) \Rightarrow (\ref{ineq2})$ yields that
\[
F_n\cap G_n \subset \{\phi(b_n,Y_n) > c(1-\delta)\psi(r_n,Z_n)\}.
\]
Therefore, by inequalities (\ref{bounds}),
\[
F_n\cap G_n 
\subset \left\{ \hat p_{\bar k}(X^n)\geq c(1-\delta)\hat q_{\bar k}(X^n) \right\}.
\]
It follows from Lemma~\ref{la3} and the definition of $d$ that for every 
$\vartheta=(p,q)$ with $p_{\bar k}\geq c q_{\bar k}$, $\prob_\vartheta (F_n^c) \leq \frac \epsilon 2$.
By Lemma~\ref{la4}, there exists $N_2\in\nats$ so that $\prob _\vartheta(G_n^c) \leq\frac\epsilon 2$
for all $n$ with $b_n p_{\bar k} \geq N_2$.
Thus, if $p_{\bar k}\geq c q_{\bar k}$, $n\in{\mathcal N}$ and $b_np_{\bar k}\geq N := \max (N_1,N_2)$, then
\[
\prob_\vartheta ( \hat p_{\bar k}(X^n)\geq c(1-\delta)\hat q_{\bar k}(X^n))
\geq 1- \prob_\vartheta(F_n^c) - \prob_\vartheta(G_n^c) \geq 1-\epsilon. 
\]
Note that $N$ does not depend on the sequence of the choices of the dice. 
$\Box$

\begin{rem}
If $K=2$, then for every $n\geq 1$ and every fixed number of times the red 
die is chosen in the first $n$ periods, the Bayes estimate of $q_{\bar k}$ 
can be shown to be an increasing function of the number of times the red 
die lands on side $\bar k$. This fact can be combined with  
Theorem~\ref{thm1} 
to give an alternative proof of Theorem~\ref{thm2} for the case $K=2$. The monotonicity 
result does not hold for $K>2$ and our proof of Theorem~\ref{thm2} does not 
use Theorem~\ref{thm1}.
\end{rem}

{\em Proof of Corollary~\ref{cor1}.}
By Chebyshev's inequality, $[B_n-\expe(B_n)]/n$ converges in probability to 0.
Thus, by condition (\ref{EBncondition}), there exists $\eta>0$ and $N_1\in\nats$
so that the event $F_n=\{\eta \leq B_n / n\leq 1-\eta\}$ has probability
$\prob(F_n)\geq 1-\frac \epsilon 2$ for all $n\geq N_1$.
By Theorem~\ref{thm2}, there exists $N_2\in\nats$ so that
\[
\prob_\vartheta (\hat p_{\bar k} (X^n) \geq c (1-\delta) \hat q_{\bar k}(X^n)|B_n=b_n) \geq 1-\frac \epsilon 2
\] 
for all $\vartheta=(p,q)\in\Delta^2$ with $p_{\bar k} \geq c q_{\bar k} $ and 
all $n\in\nats$ and $b_n\in\{1,\dots,n\}$ with $\prob(B_n=b_n)>0$ and $b_n p_{\bar k} \geq N_2$ and 
$b_n/n \leq 1-\eta $.
Let $N=\max(N_1, \lceil N_2/\eta\rceil)$. Then for every $\vartheta$ with $p_{\bar k} \geq c q_{\bar k} $ and 
every $n\in\nats$ with $n p_{\bar k} \geq N$, $F_n\subset\{B_np_{\bar k}\geq N_2, B_n/ n\leq 1-\eta\}$, so
that 
\[
\prob_\vartheta (\hat p_{\bar k} (X^n) \geq c (1-\delta) \hat q_{\bar k}(X^n)|F_n) \geq 1-\frac \epsilon 2,
\]
which implies (\ref{cor1ineq}) because $\prob(F_n) \geq 1-\frac \epsilon 2$. $\Box$

\bigskip

{\bf Acknowledgments.} 
We thank three referees for many useful suggestions. We 
thank Gary Chamberlain, Martin Cripps, Ignacio Esponda, and Muhamet Yildiz for helpful conversations.
This research is supported by National Science Foundation Grant SES 1558205.

\bigskip

\begin{center}\textbf{References}\end{center}

\medskip{}

\noindent Bochkina, N. A. and Green, P. J. (2014).
The Bernstein-von Mises theorem and nonregular models.
{\em Ann.\ Statist.} {\bf 42} 1850-1878.

\medskip{}

\noindent Diaconis, P. and Freedman, D. (1990). On the uniform consistency
of Bayes estimates for multinomial probabilities. {\em Ann.\ Statist.}
{\bf 18} 1317-1327.

\medskip{}

\noindent Dudley, R. M. and Haughton, D. (2002).
Asymptotic normality with small relative errors of posterior probabilities
of half-spaces. {\em Ann.\ Statist.} {\bf 30} 1311-1344.

\medskip{}

\noindent Fudenberg, D. and He, K. (2017). 
Type-Compatible Equilibria in Signalling Games. arXiv:1702.01819.

\medskip{}

\noindent Thompson, L. A. (2014). Bayesian Methods for Making Inferences
about Rare Diseases in Pediatric Populations. Presentation at U.S.
Food and Drug Administration. 

\medskip{}

\noindent U.S. Food and Drug Administration (2000). Guidance for the
Use of Bayesian Statistics in Medical Device Clinical Trials. Rockwell,
MD.

\medskip{}


\newpage

\begin{center}
\Large Supporting Information for
\end{center}
\begin{center}
\Large Bayesian Posteriors For Arbitrarily Rare Events
\end{center}

\vspace*{5mm}

\begin{center}
Drew Fudenberg, Kevin He, and Lorens A.\ Imhof
\end{center}

\vspace*{5mm}

\noindent{\bf Proofs of the results in Section~\ref{sectiononemultinomial}.}
Before we prove Theorem~\ref{thm1} we prove Proposition~\ref{prop1} and then state
and prove a simple consequence of Chernoff's inequality. Both results are needed 
in the proof of Theorem~\ref{thm1}.

\medskip

{\em Proof of Proposition~\ref{prop1}.} The assumption that $\pi(p)/\prod_{i=1}^K p_i^{\alpha_i-1}$ is uniformly continuous 
on $\inte\Delta$ implies that the function has a continuous extension $\tilde\pi :\Delta \to\reals$,
see Dugundji (1966), Theorem 5.2, page 302. Let $\tilde\pi_0 =\min\{\tilde\pi(p) : p\in\Delta\}$. 
Then $\tilde\pi_0>0$. Given $\epsilon>0$, choose $\delta\in(0,\tilde\pi_0)$ so small that 
\begin{equation}\label{deltachoice}
\frac{1+\frac \delta{\tilde\pi_0}}{1-\frac\delta{\tilde\pi_0}} \leq 1+\epsilon.
\end{equation}
To approximate the integrals in the assertion by sums of Dirichlet integrals we use the fact 
that the continuous function $\tilde\pi$ can be uniformly approximated by Bernstein polynomials, 
see Lorentz (1986), pages 6 and 51. Thus, there is a polynomial 
\[
h(p)=\sum_{\substack{\nu_1,\dots,\nu_K\geq 0\\ \nu_1+\dots+\nu_K=m}}c_\nu \prod_{i=1}^Kp_i^{\nu_i},
\qquad
c_\nu= \tilde\pi\left(\frac {\nu_1} m,\dots,\frac{\nu_K}m\right)\frac{m!}{\nu_1! \cdots \nu_K!},
\]
so that
\[
|\tilde \pi(p)- h(p)| \leq \delta, \qquad p \in\Delta.
\]
Using the formula 
\[
\int \prod_{i=1}^K p_i^{s_i-1}\, d\lambda(p)
=\frac{\prod_{i=1}^K \Gamma(s_i)}{\Gamma(\sum_{i=1}^Ks_i)},\qquad
s_1,\dots,s_K>0,
\]
and the relation $\Gamma(s+1)=s\Gamma(s)$, we get
\begin{align*}
\frac{\int p_k\left(\prod_{i=1}^K p_i^{s_i-1}\right)h(p)\, d\lambda(p)} 
{\int \left(\prod_{i=1}^K p_i^{s_i-1}\right)h(p)\, d\lambda(p)}
&=
\frac 1 {m+\sum_{i=1}^K s_i}
\frac{\sum_{\nu} c_\nu (\nu_k+s_k)\prod_{i=1}^K\Gamma(\nu_i+s_i)}
{\sum_\nu c_\nu \prod_{i=1}^K\Gamma(\nu_i+s_i)} .
\end{align*}
Since $c_\nu>0$ for every $\nu$, it follows that
\begin{equation}\label{bounds1}
\frac{s_k}{m + \sum_{i=1}^K s_i} 
\leq 
\frac{\int p_k\left(\prod_{i=1}^K p_i^{s_i-1}\right)h(p)\, d\lambda(p)} 
{\int  \left(\prod_{i=1}^K p_i^{s_i-1}\right)h(p)\, d\lambda(p)}
\leq 
\frac{m+s_k}{m+\sum_{i=1}^K s_i} .
\end{equation}

For all $p\in\Delta $,  $h(p)\geq \tilde\pi_0$, and so
$|\tilde \pi(p) - h(p)| \leq \delta \leq \frac{\delta}{\tilde\pi_0}h(p)$. Thus, 
\[
\left(1-\frac\delta{\tilde\pi_0}\right) h(p)
\leq \tilde\pi(p) 
\leq \left(1+\frac{\delta}{\tilde\pi_0}\right)h(p).
\]
It follows from these inequalities together with (\ref{deltachoice}) and (\ref{bounds1}) 
that for $n, n_1,\dots,n_K\in\nats_0$ with $\sum_{i=1}^Kn_i=n$, 
\begin{align*}
\frac{\int p_k \left(\prod_{i=1}^K p_i^{n_i} \right)\pi(p)\, d\lambda(p)}
{\int \left(\prod_{i=1}^K p_i^{n_i}\right) \pi(p)\, d\lambda(p)}
&= 
\frac{\int  p_k \left(\prod_{i=1}^K p_i^{n_i+\alpha_i-1}\right) \tilde\pi(p)\, d\lambda(p)}
{\int \left(\prod_{i=1}^K p_i^{n_i+\alpha_i-1}\right) \tilde\pi(p)\, d\lambda(p)}\\
&\leq
\frac {1+\frac{\delta}{\tilde\pi_0}}{1-\frac\delta{\tilde\pi_0}}
\frac{\int p_k \left(\prod_{i=1}^K p_i^{n_i+\alpha_i-1}\right) h(p)\, d\lambda(p)}
{\int \left(\prod_{i=1}^K p_i^{n_i+\alpha_i-1}\right) h(p)\, d\lambda(p)}\\
&\leq
(1+\epsilon) \frac{m+n_k+\alpha_k}{m+n+\sum_{i=1}^K\alpha_i}.
\end{align*}
Similarly, using the inequality $1/(1+\epsilon)> 1-\epsilon$, we obtain
\begin{align*}
\frac{\int p_k \left(\prod_{i=1}^K p_i^{n_i}\right) \pi(p)\, d\lambda(p)}
{\int \left(\prod_{i=1}^K p_i^{n_i}\right) \pi(p)\, d\lambda(p)}
&\geq
\frac {1-\frac{\delta}{\tilde\pi_0}}{1+\frac\delta{\tilde\pi_0}}
\frac{\int p_k \left(\prod_{i=1}^K p_i^{n_i+\alpha_i-1}\right) h(p)\, d\lambda(p)}
{\int \left(\prod_{i=1}^K p_i^{n_i+\alpha_i-1}\right) h(p)\, d\lambda(p)}\\
&\geq
(1-\epsilon) \frac{n_k+\alpha_k}{m+n+\sum_{i=1}^K\alpha_i}.
\end{align*}
The assertion follows with $\gamma=m+\sum_{i=1}^K\alpha_i$. $\Box$

\bigskip 

\noindent
\textbf{Remark $\mathbf 3'$.} Using results on the degree of approximation by Bernstein polynomials,
one may compute explicit values for the constant $\gamma$ in Proposition~\ref{prop1}.
If, for example, $K=2$ and $\phi(p_1)=\tilde\pi(p_1,1-p_1)$ has a continuous derivative
on $[0,1]$, one can apply Theorem 1.6.1 in Lorentz (1986) to show that 
(\ref{prop1ineq}) holds with 
\[
\gamma=\alpha_1+\alpha_2+ \left\lceil \frac 54 \left( 1+\frac 2\epsilon\right)
\frac{\max\{|\phi'(p_1)| : 0\leq p_1\leq 1\}} { \min\{\phi(p_1) : 0\leq p_1\leq 1\} }
\right\rceil^2. 
\]
If $K\geq 2$ and $\tilde\pi$ coincides with a polynomial on $\Delta$, then,
by a result of Handelman (1988), $\pi$ can be written as a finite mixture of densities of Dirichlet 
distributions and Theorem 3 of Powers and Reznick (2001) gives a
computable upper bound on the support of the mixing distribution.
Thus, the inequalities in (\ref{mixineq}) hold with computable constants $a$ and $A$.

\begin{lemma}\label{la1}
Let $S_n$ be a binomial random variable with parameters $n$ and $p$.
Let $1<c<2$ and $d>0$. Then
\[
\prob\left(\frac {S_n}n \geq cp +\frac dn\right) \leq e^{(1-c)d},
\qquad
\prob\left( \frac{S_n}n \leq \frac pc -\frac dn\right) \leq e^{(1-c)d}.
\]
\end{lemma}
{\em Proof.} By Chernoff's inequality,
\[
\prob\left(\frac{S_n}n\geq cp +\frac dn\right) 
\leq 
\inf _{t>0} \left[e^{-t(cp+\frac dn)}(1-p+p e^t)\right]^n
\leq
e^{(1-c)d}[\psi(p)]^n,
\]
where $\psi(s)= e^{(1-c)cs}(1-s+se^{c-1})$. For $0\leq s\leq 1$, 
\[
\frac{\psi'(s)}{e^{(1-c)cs}}
=e^{c-1}-1 -(c-1)c - s(c-1)c(e^{c-1}-1) 
\leq
e^{c-1}-1 - (c-1)c .
\]
Set $\phi(u)=e^{u-1}-1-(u-1)u$. The function $\phi'$ is convex,
$\phi'(1)=0$ and $\phi'(2)<0$. Thus, $\phi'$ is negative on $(1,2)$,
so that $\phi(c)<\phi(1)=0$. It now follows that
$\psi$ is decreasing on $[0,1]$, so that $\psi(p)\leq\psi(0)=1$.
This proves the first claim. The proof of the second claim is similar. 
$\Box$

\bigskip

{\em Proof of Theorem~\ref{thm1}.} Let $0<\epsilon<1$. Choose $c\in(1,2)$ and $\delta>0$
so that
\[
\frac {1-\delta}c > 1-\frac\epsilon 2,\qquad (1+\delta)c < 1+\frac\epsilon 2.
\]
Let $d>0$ be so that the bound in Lemma~\ref{la1} satisfies $e^{(1-c)d} < \frac \epsilon 2$.
By Proposition~\ref{prop1}, there exists $\gamma>0$ so that for every $n\in\nats$,
\[
(1-\delta)\frac{X_k^n}{n+\gamma} \leq \hat p_k (X^n) \leq (1+\delta) \frac{X_k^n+\gamma}n,\qquad
k=1,\dots,K.
\]
Let $N$ be so large that 
\[
(1-\delta) \left(\frac 1c - \frac dN \right) \frac{1}{1+\gamma/N} > 1-\epsilon,\quad
(1+\delta) \frac{d+\gamma}N < \frac \epsilon 2.
\]
Fix $k$, $p_k$ and $n$ with $n p_k \geq N$. Set $A=\{\frac 1n X_k^n < c p_k  + \frac dn\}$
and $B=\{\frac 1n X_k^n > \frac {p_k} c - \frac dn\}$.
On $A$,
\[
\frac{\hat p_k(X^n)}{p_k} \leq (1+\delta) \frac{X_k^n+\gamma}{n p_k}
\leq (1+\delta) \left(c+\frac {d+\gamma}{n p_k}  \right)
\leq (1+\delta) \left(c+\frac {d+\gamma}{N}  \right) <1+\epsilon
\]
and on $B$,
\begin{multline*}
\frac{\hat p_k(X^n)}{p_k} 
\geq
(1-\delta) \frac {X_k^n}{ n p_k} \frac {n }{n+\gamma}
\geq
(1-\delta) \left(\frac 1c -  \frac d{n p_k} \right) \frac{1}{1+\gamma/n}\\
\geq(1-\delta) \left(\frac 1c -  \frac d N \right) \frac{1}{1+\gamma/N} > 1-\epsilon. 
\end{multline*}
By Lemma~\ref{la1}, $\prob_p(A\cap B)\geq 1-\prob_p(A^c)-\prob_p(B^c)\geq 1-\epsilon$.
$\Box$

\bigskip

\noindent
{\bf Remark $\mathbf 3''$.} In the proof of Theorem~\ref{thm1} one can choose
$c=1+\frac\epsilon 4$, $\delta=\frac \epsilon 5$, and $d=3\epsilon^{-2}$.
If the prior-dependent constant $\gamma>0$ is so chosen that the inequalities in 
(\ref{prop1ineq}) hold with $\epsilon$ replaced by $\frac\epsilon 5$, then
it follows by a small variation of the above proof that the conclusion of 
Theorem~\ref{thm1} holds for $N=8\epsilon^{-3} + 3\gamma\epsilon^{-1}$.

\bigskip

The proof of Example~\ref{ex1} uses the following lower bound for the Bayes estimates of $p_1$.

\begin{lemma}\label{la2}
Let $\pi(p)=e^{-1/p}$, $0<p\leq 1$. Then 
\[
\frac{\int_0^1 p^{\nu+1}(1-p)^{n-\nu}\pi(p)\, dp}{\int_0^1 p^\nu(1-p)^{n-\nu}\pi(p)\, dp} 
\geq \frac 1{ 8\sqrt {1\vee n}}
\]
for every $n\in\nats_0$ and $\nu=0,\dots,n$. 
\end{lemma}
{\em Proof.}
Let $U$ be a random variable with density proportional to $p^\nu(1-p)^{n-\nu}\pi(p)$
and let $V$ be a random variable with density proportional to $(1-p)^n\pi(p)$, $0<p<1$.
Then $U$ is larger than $V$ in the likelihood ratio order since
$p^\nu(1-p)^{n-\nu}\pi(p) / [(1-p)^n\pi(p)] = (p/(1-p))^\nu $ is increasing in $p$.
This implies that $\expe(U)\geq \expe(V)$, that is,
\[
\frac{\int_0^1  p^{\nu+1}(1-p)^{n-\nu}\pi(p)\, dp}{\int_0^1 p^\nu(1-p)^{n-\nu}\pi(p)\, dp}
\geq
\frac{\int_0^1  p (1-p)^n\pi(p)\, dp}{\int_0^1 (1-p)^n\pi(p)\, dp},
\]
see Lehmann and Romano (2005), page 70. It is therefore enough to prove the claim for $\nu=0$.

Let $f_n(p)=c_n(1-p)^n\pi(p)$, where $c_n=[\int_0^1 (1-p)^n\pi(p)\, dp]^{-1}$. We have
\[
f_n'(p)=c_n\frac{e^{-1/p}(1-p)^{n-1}}{p^2}(1-p-np^2),
\]
showing that $f_n$ is increasing on $[0,2a_n]$, where $a_n=1/(4\sqrt{ 1\vee n})$. Hence
\[
\frac{\int_{a_n}^1 f_n(p)\, dp}{1-\int_{a_n}^1f_n(p)\, dp}
=\frac{\int_{a_n}^1 f_n(p)\, dp}{\int_0^{a_n}f_n(p)\, dp} 
\geq \frac{\int_{a_n}^{2a_n} f_n(p)\, dp}{a_n f_n(a_n)}
\geq \frac{(2a_n-a_n)f(a_n)}{a_n f_n(a_n)}
= 1.
\]
Thus $\int_{a_n}^1 f_n(p)\, dp \geq \frac 12$, and therefore
\[
\int_0^1 pf_n(p)\, dp 
\geq \int_{a_n}^1 p f_n(p)\, dp 
\geq a_n \int_{a_n}^1 f_n(p)\, dp 
\geq \frac 12 a_n
=\frac 1{8\sqrt{1\vee n}}. \ \Box
\]

\bigskip

{\em Proof of Example~\ref{ex1}.} Let $N\in\nats$. For $n>N^2$ define 
$p(n)\in\Delta$ by $p_1(n)=N n^{-\frac 12-\delta}$. By Lemma~\ref{la2}, 
$\hat p_1(X^n) -2 p_1(n) \geq  n^{-\frac 12 } (\frac 18-2 N n^{-\delta})$,
and so, for $n$ sufficiently large, $\prob_{p(n)}(|\hat p_1(X^n) - p_1(n)| > 
p_1(n))=1$. $\Box$

\bigskip

{\em Proof of Example~\ref{ex2}.} Suppose $\pi$ satisfies Condition 
${\mathcal P}(\alpha)$, $\alpha\in(0,\infty)^K$. By Proposition~\ref{prop1}, there exists 
$\gamma>0$ so that $\hat p_1(X^n) \geq \alpha_1/[2(n+\gamma)]$. For every $n>\alpha_1/8$ 
pick $p(n)\in\Delta$ with $p_1(n)=\alpha_1/(8n)$. Let $n_0=\max (\alpha_1/8, \gamma)$.
If $n> n_0$, then $\alpha_1/[2(n+\gamma)] > 2 p_1(n)$, and so
$\prob_{p(n)}(|\hat p_1(X^n)-p_1(n)| > p_1(n))=1$. 
Since $\limsup_{n\to\infty}\zeta(n)/n=\infty$, there exists for every $N\in\nats$ an $n>n_0$
with $\zeta(n)p_1(n) \geq N$. $\Box$

\bigskip
 
The following result was used in Remark~\ref{rem4}.

\begin{proposition}\label{prop2}
Let $K>2$ and $\bar k\in\{1,\dots, K\}$.
Suppose the density $\pi$ of the prior distribution on $\Delta$ 
satisfies Condition ${\mathcal P}(\alpha_1,\dots,\alpha_K)$
with $\alpha_1,\dots,\alpha_K>0$.  
Then the image measure induced by the mapping $(p_1,\dots,p_K)\mapsto
(p_{\bar k},\sum_{k\neq \bar k}p_k)$ has a density that satisfies Condition 
${\mathcal P}(\alpha_{\bar k}, \sum_{k\neq \bar k}\alpha_k)$.
\end{proposition}
{\em Proof.} Suppose without loss of generality that $\bar k=1$.
Then the image measure has a density $\pi_1$ with respect to the normalized Lebesgue measure on $\Delta_1=\{q\in[0,1]^2:q_1+q_2=1\}$
which is given by
\[
\pi_1(q) = \int _{A(q_2)} \pi\biggl(q_1,p_2,\dots,p_{K-1},q_2- \sum_{k=2}^{K-1}p_k \biggr)\, d(p_2,\dots,p_{K-1}),
\]
where 
\[
A(q_2) = \{(p_2,\dots,p_{K-1}) \in (0,1)^{K-2} : p_2+\dots+p_{K-1} < q_2\}.
\]
Making the change of variable $t=(t_2,\dots,t_{K-1}) =  q_2^{-1}(p_2,\dots,p_{K-1})$ we get
\[
\pi_1(q)=q_2^{K-2}\int_{A(1)} \pi\biggl( q_1, q_2t, q_2\biggl(1-{\sum_{k=2}^{K-1}t_k } \biggr)\biggr)\, dt
\]
for $q\in\Delta_1$ with $q_2>0$.
Since $\pi$ satisfies Condition ${\mathcal P}(\alpha_1,\dots,\alpha_K)$, 
there exists a continuous positive function $\tilde\pi$ on $\Delta$ such that
$\tilde\pi(p)=\pi(p)/\prod_{k=1}^K p_k^{\alpha_k-1}$
for all $p\in\inte\Delta$. Hence, for $q\in\inte\Delta_1$, 
\begin{multline*}
\frac{\pi_1(q)}{q_1^{\alpha_1-1}q_2^{(\sum_{k= 2}^K \alpha_k)-1}}\\
=\int_{A(1)} \tilde\pi\biggl(q_1, q_2 t, q_2\biggl(1-\sum_{k=2}^{K-1}t_k\biggr)\biggr) \prod_{k=2}^{K-1}t_k^{\alpha_k-1} 
\biggl(1-\sum_{k=2}^{K-1} t_k\biggr)^{\alpha_K-1}\, dt .
\end{multline*}
The integral is positive for every $q\in\Delta_1$ and, by dominated convergence, 
depends continuously on $q\in\Delta_1$.
Thus, $\pi_1$ satisfies condition ${\mathcal P}(\alpha_1,\alpha_2+\dots+\alpha_K)$. $\Box$

\vspace*{5mm}

\noindent{\bf Proofs of the examples in Section~\ref{sectiontwomultinomials}.}

\medskip

{\em Proof of Example~\ref{ex3}.}
Let $N\in\nats$. For every $n\geq \max( N , \frac Nc)$ let 
$p(n)=(\frac Nn, 1-\frac Nn)$, $q(n)=(\frac N{cn}, 1-\frac N{cn})$, $\vartheta_n=(p(n),q(n))$,
and 
\[
A_n=\left\{ \hat p_1(X^n)\geq \frac c2 \hat q_1(X^n)\right\}.
\]
We will prove more than is stated, namely that $\prob_{\vartheta_n}(A_n)\to 0$ as $n\to\infty$.
Let $Y_n$ denote the number of times the blue die lands on side $1$ in the first $n$ periods.
By Proposition~\ref{prop1}, there exists $\gamma>0$ so that 
$\hat p_1(X^n) \leq \frac 32 (Y_n+\gamma)/(B_n+\gamma)$.
For every $n\geq \max(N,\frac Nc)$ and $b\in\{0,1,\dots,n\}$, by Lemma~\ref{la2},
\[
\prob_{\vartheta_n}(A_n|B_n=b) 
\leq 
\prob_{\vartheta_n}\left( \left. \frac 32\frac{Y_n +\gamma}{ b+\gamma} \geq  \frac c{16\sqrt{1\vee (n-b)}}\right| B_n=b\right).
\]
If $b>\frac n2\mu_B$, then 
$c(b+\gamma)/(24\sqrt{1\vee (n-b)})\geq  d\sqrt n$ with 
$d:=c\mu_B/(48\sqrt{1-\mu_B/2})$, and it follows that
\[
\prob_{\vartheta_n}(A_n|B_n=b) \leq \prob_{\vartheta_n} (Y_n \geq -\gamma+ d \sqrt n | B_n=b).
\]
To bound the probability on the right-hand side we use a Poisson approximation to the 
conditional distribution of $Y_n$. Let $W_\nu$ be a Poisson random variable with mean 
$\nu$. Then, by Stein (1986), (43) on page 89,
\begin{align*}
\prob_{\vartheta_n}\left( Y_n \geq - \gamma + d\sqrt n | B_n=b \right) 
&\leq \prob\left( W_{b p_1(n)} \geq - \gamma + d\sqrt n \right) + p_1(n)\\
&\leq \prob\left(W_N \geq -\gamma + d\sqrt n \right)+ \frac Nn.
\end{align*}
In the second line we used the fact that 
$W_N$ is stochastically larger than $W_{b p_1(n)}$ because
$N \geq b p_1(n)$,
see Lehmann and Romano (2005), pages 67-70. Hence
\begin{align*}
\prob_{\vartheta_n}(A_n) 
& \leq \prob_{\vartheta_n} \left( B_n\leq \frac n2\mu_B\right) +\sum_{b: b>\frac n2\mu_B} 
\prob_{\vartheta_n}(A_n|B_n=b) \prob_{\vartheta_n}(B_n=b)\\
&\leq \prob_{\vartheta_n} \left(\frac 1n B_n\leq \frac 12\mu_B\right) + 
\prob\left( W_N \geq -\gamma + d \sqrt n\right)+ \frac Nn.
\end{align*}
As $n\to \infty$, $\prob(W_N \geq -\gamma + d \sqrt n )\to 0$ and, by the weak law of large 
numbers, $\prob_{\vartheta_n}(\frac 1n B_n\leq \frac 12\mu_B)\to 0$.  Thus, 
$\prob_{\vartheta_n}(A_n)\to 0$ as $n\to\infty$. $\Box$

\bigskip

{\em Proof of Example~\ref{ex4}.}
Let $Y_n$ and $Z_n$ be the respective number of times the blue and the red die land 
on side $1$ in the first $n$ periods. By Proposition~\ref{prop1}, there exists $\gamma>0$ 
so that
\begin{align*}
\prob_{\vartheta}\left( \hat p_1(X^n) < \frac c2 \hat q_1(X^n) \right) 
& \geq
\prob_{\vartheta}\left( \frac 32 \frac{Y_n+\gamma}{B_n+\gamma} < \frac c4 \frac{Z_n}{n+\gamma}\right) \\
&\geq
\prob_\vartheta\left( Y_n=0, \frac{ 6 \gamma}c  < \frac{B_n}{n}Z_n\right) .
\end{align*}
For every $n\in\nats$ with $n \geq c$ pick $\vartheta_n=(p(n),q(n))\in\Delta^2$ 
with $p_1(n)=\frac cn$ and $q_1(n)=\frac 1n$. Let $\mu_0\in(0,\mu_B)$ and 
$\mu_1 \in (\mu_B,1)$. Then, for $b=  \lceil\mu_0 n\rceil , \dots,\lfloor\mu_1 n\rfloor $,
\[
\prob_{\vartheta_n}\left( \left. Y_n=0, \frac{ 6 \gamma}{c } < \frac{B_n} n Z_n \right |B_n=b \right)
\geq 
[1-p_1(n)]^n 
\prob_{\vartheta_n} \left( \left. \frac{ 6 \gamma}{c\mu_0 } < Z_n \right| B_n=\lfloor \mu_1 n \rfloor \right).
\]
Now $[1-p_1(n)]^n\to e^{-c}>0$ and, by Stein (1986), (43) on page 89,
\[
\prob_{\vartheta_n}\left( \left. \frac{ 6 \gamma}{c \mu_0} < Z_n \right| B_n=\lfloor \mu_1 n \rfloor \right) 
\geq 
\prob \left(W > \frac {6\gamma}{c\mu_0 } \right)-\frac 1n,
\]
where $W$ is a Poisson random variable with mean $1-\mu_1$. Hence 
\[
\liminf _{n\to\infty} 
\prob_{\vartheta_n}\left( \left. Y_n=0, \frac{6\gamma}c < \frac{B_n} n Z_n \right| 
\mu_0 n \leq B_n \leq \mu_1 n\right) 
>0.
\]
Since $\prob(\mu_0 n \leq B_n \leq \mu_1 n)\to 1$, it follows that there exists 
$\epsilon_0>0$ and $n_0\in\nats$ so that
\[
\prob_{\vartheta_n}\left( \hat p_1(X^n) < \frac c2 \hat q_1(X^n) \right) 
> \epsilon_0
\]
for all $n\geq n_0$. Since $\zeta(p_1(n))/p_1(n)\to\infty$
as $n\to\infty$, there exists for every $N\in\nats$ an $n\geq n_0$ with 
$n\zeta(p_1(n)) \geq N$ and $\vartheta_n$ has the required properties. $\Box$

\vspace*{5mm}

\noindent{\bf Proofs of the auxiliary results in Section~\ref{proofofthm2}.}

\medskip

{\em Proof of Lemma~\ref{la3}.}
Set $\ell=d/(n\wedge m)$. By Markov's inequality, for every $t>0$,
\begin{equation}\label{la3eq1}
\prob\left(\frac{T_m}m \geq \frac 1{c'} \frac{S_n} n + \ell \right)
=
\prob\left( e^{t(c'T_m - \frac m n S_n)} \geq e^{t c' \ell m}\right)
\leq
\frac{\expe[e^{t(c'T_m - \frac m n S_n)} ]}{e^{t c' \ell m}}.
\end{equation}
We will determine a suitable value for $t$ so that the expectation is at most $1$. 
Let $\xi$ and $\tau$ be Bernoulli variables with $\prob(\xi=1)=p$ and $\prob(\tau=1)=q$. Then
\begin{equation}\label{la3eq2}
\expe[e^{t(c'T_m - \frac m n S_n)}] =
\expe(e^{t c' T_m})  \expe (e^{- t \frac m n S_n}) =
[\expe(e^{t c' \tau})]^m [\expe( e^{- t \frac m n \xi})]^n.
\end{equation}
For $t>0$ and $s\in\reals$ let $\psi_t(s)=(1-s+s e^{c't}) (1-cs + cs e^{-t})$. Since $p\geq c q$,
\[
\expe (e^{tc' \tau}) \expe(e^{-t\xi}) = (1-q+q e^{c't}) (1-p + p e^{-t})
\leq \psi_t(q).
\]
We have $\psi_t(0)=1$, and $\psi_t''(s)=2c(e^{c't}-1)(e^{-t}-1)<0$, so that
$\psi_t$ is concave. For $t_0:=(c'+1)^{-1}\log (c/c')$,
\[
\psi_{t_0}'(0)=e^{c't_0}-1+ c( e^{-t_0}-1) = \int_0^{t_0} e^{-u} [c' e^{(c'+1)u} - c]\, du <0,
\]
so that $\psi_{t_0}(s)\leq 1$ for $s\geq 0$. Hence,
\begin{equation}\label{la3eq3}
\expe(e^{c't_0\tau}) \expe(e^{-t_0\xi}) \leq 1. 
\end{equation}

If $m\leq n$, then by Lyapunov's inequality, 
$[\expe(e^{- t_0 \frac m n \xi })]^n \leq  [ \expe(e^{- t_0 \xi}) ]^m$. 
Combining this inequality with (\ref{la3eq2}) and (\ref{la3eq3}) yields
\[
\expe[e^{t_0(c' T_m  - \frac m n S_n)}] \leq 
[\expe(e^{t_0 c'\tau  })]^m [\expe(e^{- t_0 \xi})]^m
\leq  1, 
\]
and so, by (\ref{la3eq1}), 
\[
\prob\left(\frac{T_m}m \geq \frac 1{c'}\frac{S_n} n + \ell \right) \leq {e^{-t_0 c' \ell m}} 
= \left( \frac {c'}c \right)^{c'd /(c'+1)} .
\]

If $m > n$, then Lyapunov's inequality gives 
$[\expe(e^{t c'\tau })]^m \leq [\expe (e^{t c' \frac  m n\tau})]^n$.
Setting $t_1= \frac n m t_0$, we get in this case
\[
\expe[e^{t_1( c' T_m - \frac m n S_n)}] 
\leq
[ \expe( e^{t_1 c' \frac m n \tau})]^n [\expe(e^{- t_1 \frac m n \xi})]^n
\leq 1,
\]
and so 
\[
\prob\left(\frac{T_m}m \geq \frac 1{c'} \frac{S_n}n + \ell \right) \leq e^{- t_1 c' \ell m} 
= \left(\frac {c'}c\right)^{\displaystyle c'd/(c' + 1)}. \ \Box
\]

\medskip

{\em  Proof of Lemma~\ref{la4}.} We will use a Poisson approximation to the binomial distribution.
If $W_\nu$ is a Poisson random variable with mean $\nu>0$, then
$\prob(W_\nu\leq M) \to 0$ 
as $\nu\to \infty$. Thus there exists $N_0\in\nats$ so that 
$\prob(W_\nu \leq M) <  \frac 12\epsilon$ for $\nu>N_0$.
By Stein (1986), (43) on page 89, $|\prob_p(S_n \leq M) - \prob(W_{np}  \leq M)|\leq p$. 
Thus if $np\geq N_0$ and $p\leq\frac 12 \epsilon$, then $\prob_p(S_n \leq M) \leq \epsilon$.
In particular, for $p=\frac 12\epsilon$ and $n=\lceil 2 N_0/\epsilon\rceil$, we have
$\prob_{\epsilon/2} (S_{\lceil 2 N_0/\epsilon\rceil } \leq M) \leq \epsilon$.

On the other hand, if $p > \frac 12\epsilon$ and $n\geq 2N_0/\epsilon$, then
\[
\prob_p(S_n \leq M) \leq \prob_{\epsilon/2}(S_n\leq M) 
\leq \prob_{\epsilon/2} (S_{\lceil 2 N_0/\epsilon\rceil } \leq M) \leq \epsilon,
\]
where we used the fact that the family of binomial distributions is stochastically increasing
in both parameters, see e.g.\ Lehmann and Romano (2005), pages 67-70. 
The claim follows with $N=2N_0/\epsilon$. $\Box$

\bigskip

\begin{center}\textbf{References}\end{center}

\medskip{}

\noindent Dugundji, J. (1966). {\em Topology.} Allyn and Bacon, Boston, MA.

\medskip

\noindent Handelman, D. (1988).
Representing polynomials by positive linear functions on compact convex polyhedra.
{\em Pacific J. Math.} {\bf 132} 35-62. 

\medskip

\noindent Lehmann, E. L. and Romano, J. P. (2005). 
{\em Testing Statistical Hypotheses}, third ed. Springer, New York.

\medskip

\noindent Lorentz, G. G. (1986). {\em Bernstein Polynomials}, second ed.
Chelsea, New York.

\medskip{}

\noindent
Powers, V. and Reznick, B. (2001). A new bound for P\'olya's theorem with applications 
to polynomials positive on polyhedra. {\em J. Pure Appl.\ Algebra} {\bf 164} 221-229. 

\medskip

\noindent Stein, C. (1986). {\em Approximate Computation of Expectations.}
Institute of Mathematical Statistics, Hayward, CA.
\end{document}